\documentclass{gtart_h}
\def\ifplaintex{\expandafter\ifx\csname documentclass\endcsname\relax}

\def\gtp{{\mathsurround=0pt\it $\cal G\mskip-2mu$eometry \&\ 
$\cal T\!\!$opology $\cal P\!$ublications}}  % GT publications

\def\recd{{\small Received:\qua\receiveddate\ifx\reviseddate\relax
\else\qquad Revised:\qua\reviseddate\fi\par}} 

\def\lognumber#1{\def\thelognumber{#1}}
\def\volumenumber#1{\def\thevolumenumber{#1}}
\def\volumeyear#1{\def\thevolumeyear{#1}}
\def\papernumber#1{\def\thepapernumber{#1}}
\def\pagenumbers#1#2{\def\startpage{#1}\def\finishpage{#2}}
\def\published#1{\def\publishdate{#1}}

\def\received#1{\def\receiveddate{#1}}
\def\revised#1{\def\reviseddate{#1}}
\def\accepted#1{\def\accepteddate{#1}}

\long\def\asciiabstract#1{\long\def\theasciiabstract{#1}}

\let\\\par\let\thelognumber\relax\let\thevolumenumber\relax
\let\thepapernumber\relax\let\thevolumeyear\relax\let\startpage\relax
\let\finishpage\relax\let\publishdate\relax\let\receiveddate\relax
\let\reviseddate\relax\let\accepteddate\relax\let\theasciititle\relax
\let\theasciiauthors\relax
\let\theasciiabstract\relax

\let\theasciiemail\relax

\ifplaintex
\font\logobig=cmssbx10 scaled 3836
\font\logomed=cmssbx10 scaled 2557
\else
\font\logobig=cmssbx10 scaled 4200
\font\logomed=cmssbx10 scaled 2800
\fi

\long\def\makeagttitle{   %%% start of definition of \makeagttitle
\count0=\startpage
\agt\hfill      %   Journal title (top left) 
\hbox to 45truept{\vbox to 0pt{\vglue -13truept{\logomed A\kern -.37em{\logobig 
T}\kern -.38em G}\vss}\hss}
\break
{\small Volume \thevolumenumber\ (\thevolumeyear)
\startpage--\finishpage\nl
Published: \publishdate}

\vglue .25truein

{\parskip=0pt\leftskip 0pt plus
1fil\def\\{\par\smallskip}{\Large\bf\thetitle}\par\medskip} \vglue
0.05truein

{\parskip=0pt\leftskip 0pt plus 1fil\def\\{\par}{\sc\theauthors}
\par\medskip}%
 
\vglue 0.03truein 

{\small\leftskip 25truept\rightskip 25truept{\bf Abstract}\stdspace\theabstract

{\bf AMS Classification}\stdspace\theprimaryclass
\ifx\thesecondaryclass\relax\else; \thesecondaryclass\fi\par
{\bf Keywords}\stdspace \thekeywords\par}\vglue 7truept

}   %%%% end of definition of \makeagttitle

\ifplaintex
\font\phead=cmsl9 scaled 950
\font\pnum=cmbx10 scaled 913
\font\pfoot=cmsl9 scaled 950
\headline{\vbox to 0pt{\vskip -4.5mm\line{\small\phead\ifnum
\count0=\startpage ISSN 1472-2739 (on-line) 1472-2747 (printed)
\hfill {\pnum\folio}\else\ifodd\count0\def\\{ }% 
\ifx\theshorttitle\relax\thetitle\else\theshorttitle\fi\hfill{\pnum\folio}
\else\def\\{ and }{\pnum\folio}\hfill\ifx\theshortauthors\relax\theauthors
\else\theshortauthors\fi\fi\fi}\vss}}
\footline{\vbox to 0pt{\vglue 0mm\line{\small\pfoot\ifnum\count0=\startpage
\copyright\ \gtp\hfill\else
\agt, Volume \thevolumenumber\ (\thevolumeyear)\hfill\fi}\vss}}
\else
\font=cmsl9 scaled 1050
\font\lnum=cmbx10 
\font=cmsl9 scaled 1050
\makeatletter
\def\@oddhead{{\small\lhead\ifnum\count0=\startpage ISSN 1472-2739 
(on-line) 1472-2747 (printed)\hfill {\lnum\number\count0}\else\ifodd\count0
\def\\{ }\ifx\theshorttitle\relax \thetitle \else\theshorttitle\fi\hfill
{\lnum\number\count0}\else\def\\{ and }{\lnum\number\count0}
\hfill\ifx\theshortauthors\relax 
\theauthors\else\theshortauthors\fi\fi\fi}}\def\@evenhead{\@oddhead}
\def\@oddfoot{\small\lfoot\ifnum\count0=\startpage\copyright\ \gtp\hfill\else
\agt, Volume \thevolumenumber\ (\thevolumeyear)\hfill\fi}
\def\@evenfoot{\@oddfoot}
\makeatother
\fi
\let\maketitlepage\makeagttitle
\let\makeshorttitle\maketitlepage
\let\maketitle\maketitlepage

\newwrite\gtoutfile
\long\gdef\makeheadfile{  %%% start of definition of \makeheadfile
{\def\\{, }\def\s{ }
\immediate\openout\gtoutfile head.xxx
\immediate\write\gtoutfile{Proxy-for: \ifx\theasciiauthors\relax
\theauthors\else\theasciiauthors\fi\s<\ifx\theasciiemail\relax\theemail\else\theasciiemail\fi>}
\immediate\write\gtoutfile{\noexpand\\}
\immediate\write\gtoutfile{Authors: \ifx\theasciiauthors\relax
\theauthors\else\theasciiauthors\fi}
{\def\\{ }\immediate\write\gtoutfile{Title: \ifx\theasciititle\relax
\thetitle\else\theasciititle\fi}}
\immediate\write\gtoutfile{Subj-class: GT or SG, GR etc}
\immediate\write\gtoutfile{MSC-class: \theprimaryclass\ifx\thesecondaryclass\relax\else, \thesecondaryclass\fi}
\immediate\write\gtoutfile{Journal-ref: Algebr. Geom. Topol. \thevolumenumber\s
(\thevolumeyear) \startpage-\finishpage}
\immediate\write\gtoutfile{Comments: Published by Algebraic and
Geometric Topology at}
\immediate\write\gtoutfile{\s\s\s  http://www.maths.warwick.ac.uk/agt/AGTVol\thevolumenumber/agt-\thevolumenumber-\thepapernumber.abs.html}
\immediate\write\gtoutfile{\noexpand\\}
\immediate\write\gtoutfile{}
\ifx\theasciiabstract\relax
\immediate\write\gtoutfile{\theabstract}\else
\immediate\write\gtoutfile{\theasciiabstract}\fi
\immediate\write\gtoutfile{}
\immediate\write\gtoutfile{\noexpand\\}
\immediate\write\gtoutfile{}
\immediate\closeout\gtoutfile}}  %%% end of definition of \makeheadfile

\def\maketitlepage{\makeagttitle\makeheadfile}
\let\makeshorttitle\maketitlepage
\let\maketitle\maketitlepage

\lognumber{47}
\volumenumber{5}
\volumeyear{2005}
\papernumber{47}
\pagenumbers{1173}{1195}
\received{24 July 2004} 
\revised{15 June 2005}
\accepted{26 July 2005}
\published{18 September 2005}

\usepackage{amssymb, amsmath, graphicx}

\def\figref#1{\hyperlink{#1anchor}{Figure~\ref*{#1}}}
\def\anchor#1{\noindent\hypertarget{#1anchor}{\smash{$\phantom{99}$}}\newline}

\newtheorem{thm}{Theorem}[section]   
\newtheorem{lem}[thm]{Lemma}         
\newtheorem{cor}[thm]{Corollary}

\theoremstyle{definition}
\newtheorem{defn}[thm]{Definition}   
\newtheorem*{rem}{Remark}            
\newtheorem*{quest}{Question}

\newcommand{\cc}{\ensuremath{\mathbb{C}}}
\newcommand{\rr}{\ensuremath{\mathbb{R}}}
\newcommand{\zz}{\ensuremath{\mathbb{Z}}}

\begin{document}

\title{Overtwisted open books from sobering arcs}                    
\authors{Noah Goodman}                  
\address{Massachusetts Institute of Technology\\Cambridge, MA 02139, USA} 
\email{ngoodman@math.utexas.edu}

\begin{abstract}  
We study open books on three manifolds which are compatible with an
overtwisted contact structure. We show that the existence of certain
arcs, called sobering arcs, is a sufficient condition for an open book
to be overtwisted, and is necessary up to stabilization by positive
Hopf-bands. Using these techniques we prove that some open books
arising as the boundary of symplectic configurations are overtwisted,
answering a question of Gay in \cite{gay}.
\end{abstract}

\asciiabstract{We study open books on three manifolds which are compatible with an
overtwisted contact structure. We show that the existence of certain
arcs, called sobering arcs, is a sufficient condition for an open book
to be overtwisted, and is necessary up to stabilization by positive
Hopf-bands. Using these techniques we prove that some open books
arising as the boundary of symplectic configurations are overtwisted,
answering a question of Gay in Algebr. Geom. Topol. 3 (2003) 569--586.}

%  AMS classification numbers, primary and secondary, and keywords :


%                      End of header form
%                      ==================
%
%   For G&T articles leave \maketitlepage uncommented.
%   For AGT and GTM articles comment it out and uncomment \makeshorttitle
%
%\maketitlepage    %%% Makes a title page for G&T articles
%
 \makeshorttitle  %%% Makes a short header for AGT and GTM articles

%
%%%%%%%%%%%%%%%%%%%%   Start of main body of article

%\maketitle

\section{Introduction}

A contact structure on a three-manifold is a nowhere integrable 
two-plane field (we consider only oriented, co-orientable, positive contact structures). 
A fundamental theorem of Eliashberg \cite{yasha} shows that contact 
structures fall into two classes: {\it overtwisted}, and {\it tight}. 
The overtwisted contact structures are determined by their homotopy 
type as plane fields, while the tight ones are more subtle.

An open book is a pair $( \Sigma , \phi )$, where $\Sigma$ is a surface 
with boundary, and $\phi$, the monodromy, is an automorphism of 
$\Sigma$. The mapping torus $ \Sigma \times_\phi S^1$ has torus boundary 
with preferred meridian and longitude, we can thus form a three manifold 
$M_{(\Sigma , \phi )}$ by gluing solid tori to the boundary. We say 
that $( \Sigma , \phi )$ is an open book on $M_{(\Sigma , \phi )}$, each 
$\Sigma \subset M_{(\Sigma , \phi )}$ is a page, while the core of 
the solid tori are the binding.

A contact structure and an open book on the same manifold are called 
{\it compatible} if the contact planes can be isotoped arbitrarily 
close to the tangent planes of the pages, while the binding remains
transverse. For example the standard tight contact structure on
$S^3$, written $\xi_0$, is compatible with the positive Hopf-band, $H^+$,
whose binding is two fibers of the Hopf fibration. 
 
 If we have two contact structures compatible with the same 
open book it is not hard to show (by first pushing the contact planes 
close to the pages, then using a straight-line homotopy) that they must 
be isotopic. (We will write $\xi_B$ for a contact structure compatible with an 
open
book $B$, unique up to isotopy.) Because of this uniqueness the contact structures 
compatible with a given open book are either all tight or 
all overtwisted.

\begin{defn}
An open book $(\Sigma , \phi )$ is called {\it overtwisted} if there is an overtwisted 
contact structure $\xi_{(\Sigma , \phi )}$ compatible with it.
\end{defn}

Giroux has shown \cite{giroux} that there is an open book compatible 
with any contact structure, so it is natural to study contact structures 
by studying open books. The first question one would like to answer is: 
when is an open book overtwisted?

A few open books can be shown to be overtwisted by general 
considerations or special techniques. 
For example, every contact structure on the Poincare homology-sphere with reversed 
orientation is overtwisted \cite{notight}, so any compatible open book must also 
be overtwisted. 
This manifold arises as $+1$-surgery on the right-handed trefoil $T^+$, so it has 
an open book which comes 
from the natural open book structure of $T^+$ by adding 
a boundary parallel negative Dehn twist.
In section \ref{sobarcs} we'll show more directly that this open book is 
overtwisted.

In this paper we study arcs (which we call {\it sobering arcs}) which 
intersect their monodromy images in specific ways.  We show that the existence of a sobering 
arc is a sufficient condition for an open book to be overtwisted. This 
is not 
a necessary condition (there are overtwisted open books with no sobering 
arc), but it is necessary up
to stabilization (by Murasugi sum with positive Hopf-bands). The most 
satisfying statement of our result in this paper is:

\begin{thm}
An open book is overtwisted if and only if it is stably equivalent to an 
open book with a sobering arc.
\end{thm}

This technique gives some interesting examples: we look at the open 
books 
that arise from positive surface configurations in symplectic manifolds, studied by Gay in \cite{gay}, and show 
that some of these open books are overtwisted, and hence the corresponding symplectic configuration cannot 
embed in a closed symplectic 4-manifold (answering a question 
of \cite{gay}).

In section \ref{background} we review some necessary background material 
on contact structures and open books. In section \ref{sobarcs} we 
introduce sobering arcs, prove that an open book with a sobering arc 
is overtwisted, and give some examples. In section \ref{specialcase} we 
investigate 
Hopf-bands and open books which have a Hopf-band as a Murasugi summand. 
In section \ref{criterion} we show that every overtwisted open book is 
stably-equivalent to one with a sobering arc. In section \ref{sympl} we 
look 
at contact manifolds that arise as the boundaries of symplectic 
configurations and show that some of these are overtwisted.

This results of this paper were contained in the author's 
Ph.D. dissertation \cite{diss}.

We would like to thank John Etnyre for a great deal of encouragement and many helpful conversations, and the anonymous reviewer for very useful comments.

\section{Background}
\label{background}

Throughout this section we will omit proofs that can be easily found in the literature. In 
particular, there are now many good references on contact geometry, see \cite{johnintro}.

\subsection{Topological conventions}

Unless we explicitly state otherwise all manifolds are three dimensional and oriented, and all 
surfaces are also oriented.

If $S \subset M$ is a surface with boundary embedded in three-manifold $M$, and $T$ another surface 
so that $\partial S \subset T$, then the {\it framing difference} for $\partial S$ between $S$ and 
$T$, written $Fr(\partial S;S,T)$ or just $Fr(S,T)$, is the oriented number of intersections between 
$S$ and a push-off of $\partial S$ along $T$. 

An automorphism of a surface is a bijective self-homeomorphism, fixing 
the boundary point-wise.
The {\it mapping torus} of a surface automorphism, $\phi : \Sigma \rightarrow \Sigma$, is the 
three-manifold $ \Sigma \times_{\phi} S^1 = (\Sigma \times I) / (p, 1) \sim (\phi (p), 0)$.

Given a simple closed curve $c \subset \Sigma$ we can define an automorphism $D^+_c : \Sigma \rightarrow \Sigma$, 
which has support only near $c$, as follows. Let $N$ be a neighborhood of $c$ 
which is identified (by oriented coordinate charts) with the annulus 
$\{ a \in \cc | 1 \leq ||a|| \leq 2 \} $ in $\cc$. Then $D_c$ is the map $a \mapsto e^{-i 2\pi (||a|| -1)}a$ 
on $N$, and the identity on $\Sigma \setminus N$. We call $D^+_c$ a positive
Dehn twist; the inverse, written $D^-_c$ is a negative Dehn twist.

\subsection{Open books}

Open books $(\Sigma_i, \phi_i )$, $i=1,2$, are {\it equivalent} if $\Sigma_1 = \Sigma_2$ and $\phi_1$ is isotopic to 
$h^{-1} \circ \phi_2 \circ h$ for some automorphism $h$.

Associated to each open book there is a three-manifold, 
$M_{\Sigma, \phi } = (\Sigma \times_{\phi} S^1) \cup_{\partial} (D^2 \times \partial  \Sigma)$, 
where the union is taken by gluing boundary tori in such a way that the curve 
$\{ p \} \times S^1$, for $p \in \partial \Sigma$, glues to meridian $\partial (D^2 \times \{ p \})$, and 
$\partial \Sigma \times \{ \theta \}$ glues to $\{ e^{i\theta} \} \times \partial \Sigma$. 
The link $\{ 0 \} \times \partial \Sigma \subset D^2 \times \partial \Sigma \subset M_{\Sigma, \phi }$ is called
the {\it binding}, sometimes we'll write it $\partial \Sigma \subset M_{\Sigma, \phi }$.
As expected, equivalent open books give homeomorphic $( M_{\Sigma, \phi }, \partial \Sigma )$.

Given two open books, and some attaching data, we can get a new open book via the Murasugi sum. Unless 
stated otherwise our Murasugi sum will be the simplest kind, defined below in terms of rectangles. A more general 
Murasugi sum, also called generalized plumbing or $*$-product in the literature, can be defined 
similarly by using polygons with $2m$ sides.

\begin{defn}
Let $(\Sigma_i , \phi_i )$, $i=1,2$, be open books, $l_i \subset \Sigma_i$ properly embedded arcs, 
and $R_i$ 
rectangular neighborhoods of $l_i$ such that $R_i \cap \partial \Sigma_i$ is two arcs of $\partial 
R_i$. The surface 
$\Sigma_1 * \Sigma_2$ is the union $\Sigma_1 \cup_{R_i} \Sigma_2$ with $R_1$ identified to $R_2$ in 
such a way 
that $\partial R_1 \cap \partial \Sigma_1 = \partial R_2 \setminus \partial \Sigma_2$, and vice versa. (This is the 
identification which gives a new surface, rather than a branched surface.)
The {\it Murasugi sum} of the $(\Sigma_i , \phi_i )$ along $l_i$ is the open book $(\Sigma_1 , \phi_1 
) * (\Sigma_2 , \phi_2 ) = (\Sigma_1 * \Sigma_2 , \phi_1 \circ \phi_2)$ (where each $\phi_i$ is 
understood to extend to $\Sigma_1 * \Sigma_2$ by the identity on $(\Sigma_1 * \Sigma_2) 
\setminus \Sigma_i$).
We will often suppress mention of the attaching arcs $l_i$ when they are clear or irrelevant.
\end{defn}

\begin{thm}
\label{muriscon}
The manifold $M_{(\Sigma_1 , \phi_1 ) * (\Sigma_2 , \phi_2 )}$ is $M_{\Sigma_1, 
\phi_1 } \# M_{\Sigma_2, \phi_2 }$.
\end{thm}
\begin{proof} See \cite{gabai}. \end{proof}

The {\it Hopf-bands}, $H^\pm$, are the open books with surface an 
annulus and monodromy a single positive or negative Dehn twist about the 
center of the annulus.

\subsection{Contact structures}

A one dimensional submanifold $L \subset (M, \xi )$ is {\it Legendrian} if $T_pL \subset \xi$ at each 
$p \in L$, and is {\it transverse} if $T_pL \not\subset \xi$. Any curve in a contact manifold can be 
made Legendrian or transverse by a $C^0$-small perturbation. 

Given an embedded surface $S \subset (M, \xi )$ the contact planes restrict to a (singular) line field 
on $S$. This singular line field integrates into a singular foliation (since any line field is 
integrable), which we call the 
{\it characteristic foliation}, $\xi |_S$.

\label{tw}
For $S \subset M^3$ a surface with $\partial S$ Legendrian, $tw(\partial 
S, S)$ is the difference in framing of $\partial S$ between $\xi$ and $S$. 
(That is the intersection number of $S$ with a push-off of $\partial S$ 
along $\xi$.)

A contact structure is {\it overtwisted} if there is an embedded disk $D$ with 
Legendrian boundary such that 
$tw(\partial D, D) = 0$.
A contact structure that is not overtwisted is {\it tight}. 

We now state some key theorems about overtwisted contact structures:

\begin{thm}[The Bennequin Inequality]
\label{bennineq}
If $(M,\xi )$ is a tight contact manifold, $S \subset M$ is an 
embedded surface, and $\partial S$ is Legendrian, then:
\[ tw(\partial S , S) \leq -\chi (S) - |r| \leq -\chi (S). \]
\end{thm}

In the previous theorem, $r$ is the {\it rotation number} of the contact 
structure over the surface. We will use only the courser inequality in this 
paper.

\begin{thm}[Eliashberg \cite{yasha}]
\label{otclass}
Overtwisted contact structures are isotopic if, and only if, they are homotopic 
as plane fields.
\end{thm}

The following theorem is due to Lutz and Martinet (see \cite{yasha}
for a review):
\begin{thm}
\label{lutzex}
There is an overtwisted contact structure in each homotopy class of plane fields.
\end{thm}

A surface $S \subset (M, \xi )$ is {\it convex} if there is a vector 
field transverse to $S$ whose flow preserves $\xi$. 
If $v$ is such a vector field, the curves $ \Gamma = \{ p \in S | v(p) \in \xi \}$ are called the {\it 
dividing curves} (or {\it divides}) of $S$. For a convex surface with Legendrian boundary $tw(\partial S, S)$ 
is negative half the number of intersection points in $\Gamma \cap \partial S$.

In \cite{girouxconvex} Giroux gives a useful characterization of convex surfaces. A characteristic 
foliation with isolated singularities is called {\it Morse-Smale} if all singularities are hyperbolic 
(in the dynamical systems sense) and there are no saddle-to-saddle connections. (There is another 
technical requirement in the definition of Morse-Smale, it will be satisfied if there are no closed 
orbits, and this will be enough for us.) A surface $S \subset (M, \xi )$ is convex 
if and only if $\xi |_S$ is Morse-Smale. For details see \cite{girouxconvex}.

\begin{lem}
Any closed surface $S \subset (M, \xi )$ can be made convex by a $C^0$-small perturbation. A surface 
with Legendrian 
boundary can be made convex by a $C^0$-small perturbation, relative to $\partial S$, if $tw(\partial 
S, S) \leq -1$ for each component of $\partial S$.
\end{lem}

We'll frequently need the following lemma, called the Legendrian Realization Principal (LeRP), due to 
Ko Honda \cite{ko}:

\begin{lem}[LeRP]
\label{lerp}
Let $c \subset S$ be a set of curves and arcs (possibly with $\partial c \subset 
\partial S$), on a convex surface $S$. If each 
component of $S \setminus c$ intersects the dividing curves of $S$, then $S$ may be perturbed, by a 
$C^0$-small perturbation through 
convex surfaces, to make $c$ Legendrian.
\end{lem}

We will often use ``LeRP" as a verb, meaning ``to apply lemma \ref{lerp} to a set of curves".

The folding trick, implicit in \cite{ko}:

\begin{lem}
\label{folding}
Let $c \subset S$ be a non-separating collection of simple closed curves on 
convex surface $S$ which don't 
intersect the divides $\Gamma \subset S$, and $c'$ a parallel copy of $c$. Then $S$ can be perturbed 
to a convex surface divided 
by $\Gamma \cup c \cup c'$. This perturbation is $C^0$-small, but not through convex surfaces. 
\end{lem}

\subsection{Compatibility}

As mentioned above, an open book is compatible with a contact structure, $\xi$ if there is a homotopy 
$\xi_t$ of plane fields (smooth for $t>0$) so that $\xi_1 = \xi$, $\xi_0$ is tangent to the pages, $\xi_t$ is contact for $t >0$, 
and the binding is transverse to $\xi_t$. 

%Ex H^\pm. Th-wink.

Torisu shows in \cite{torisu} that Murasugi sum corresponds to contact connect sum for
the compatible contact structure, that is $\xi_{B * C} = \xi_B \# \xi_C$ (where equality means contact isotopy). 
 Since $H^+$ is compatible with
the standard structure $\xi_0$, this implies (see \cite{vincent}) that $\xi_{B * H^+} = \xi_B$.  
We say that $B * H^+$ is a {\it stabilization} of $B$ and that two open books related by a series of
stabilizations are {\it stably-equivalent}.

A key theorem of Giroux relates compatible open books (\cite{giroux}, 
see \cite{diss} for an exposition of the proof):

\begin{thm}[Giroux]
\label{stabeq}
Two open books compatible with a fixed contact structure are stably-equivalent.
\end{thm}

\section{Sobering arcs}
\label{sobarcs}

\begin{figure}[ht!]\anchor{fig:intsigns}
\centering
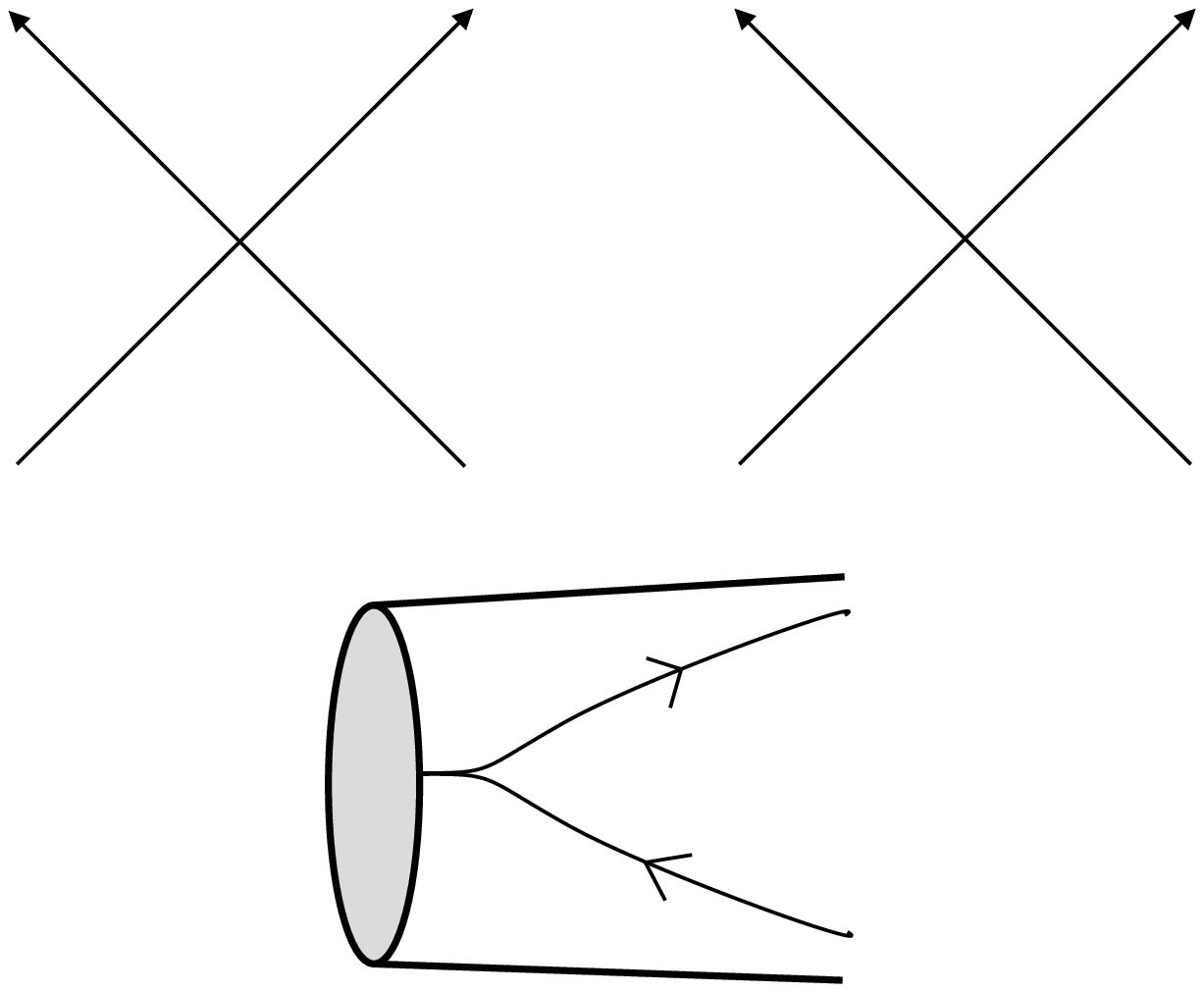
\caption{The sign conventions for the intersection 
of curves $( \alpha , \beta )$: a) negative
intersection, b) positive, c) a positive intersection as it appears at the boundary of the
surface}
\label{fig:intsigns}
\end{figure}

\begin{defn}
Let $\alpha , \beta \subset \Sigma$ be properly embedded oriented arcs 
which intersect transversely, 
$\Sigma$ an oriented surface. If $p \in \alpha \cap \beta$ and 
$T_{\alpha}$, $T_{\beta}$ are 
tangent vectors to $\alpha$ and $\beta$ (respectively) at $p$, then the 
intersection number $i_p = +1$ if the ordered basis 
$(T_{\alpha}, T_{\beta})$ agrees with the 
orientation of $\Sigma$ at $p$, and $i_p = -1$ otherwise. Let $\tilde{\beta}$ be an arc which 
minimizes intersections with $\alpha$ over boundary fixing isotopies of $\beta$.
Sign conventions are illustrated in \figref{fig:intsigns}.
We define the following intersection numbers: 
\begin{enumerate}
\item The {\it algebraic intersection number}, $i_{alg}(\alpha, 
\beta) = \sum_{p \in \alpha \cap \beta \cap int(\Sigma)} i_p$, is the oriented sum over intersections on the interior of $\Sigma$.
\item The {\it geometric intersection number}, $i_{geom}(\alpha, 
\beta) = \sum_{p \in \alpha \cap \tilde{\beta} \cap int(\Sigma)} |i_p|$, is the unsigned
count of interior intersections, minimized over all boundary fixing isotopies.
\item The {\it boundary intersection number}, $i_{\partial}(\alpha, 
\beta) =\frac{1}{2} \sum_{p \in \alpha \cap \tilde{\beta} \cap \partial\Sigma} i_p$, is one-half the oriented 
 sum over intersections at the boundaries of the arcs, after the arcs 
have been isotoped, fixing boundary, to minimize geometric intersection.
\end{enumerate}
\end{defn}

In particular, for an arc $\alpha \subset \Sigma$ in the page of an open book $(\Sigma , \phi )$ we 
may consider the intersection numbers: $i_{alg}(\alpha , \phi (\alpha ))$, $i_{geom}(\alpha , 
\phi (\alpha ))$, $i_{\partial}(\alpha , \phi (\alpha ))$. We orient $\phi (\alpha )$ 
by reversing a pushed 
forward orientation on $\alpha$. Since reversing the orientation of $\alpha$ will also reverse that 
on $\phi (\alpha )$, without changing the intersection signs, it is irrelevant 
which orientation we 
choose.

Much subtle information about a surface automorphism can be gleaned by comparing these intersection 
numbers for properly embedded arcs. For instance the algebraic monodromy (the induced map 
$\phi_* : H_*(\Sigma ) \mapsto H_*(\Sigma )$,) is unable to distinguish between $(\Sigma , \phi )$ and 
$(\Sigma , \phi^{-1})$, since $\phi_*$ is related to $\phi^{-1}_*$ by conjugation (through the 
change of basis which reverses the orientation of each basis vector). However, 
$i_{\partial}(\alpha , \phi (\alpha )) = - i_{\partial}(\alpha , \phi^{-1} (\alpha ))$ is a 
geometrically well-defined quantity which can make this distinction. 

\begin{defn}
A properly embedded arc $\alpha \subset \Sigma$ is {\it sobering}, for a monodromy $\phi$, if 
$i_{alg}(\alpha , \phi (\alpha )) + i_{\partial}(\alpha , \phi (\alpha )) + i_{geom}(\alpha , \phi 
(\alpha )) \leq 0$, and $\alpha$ is not isotopic to $\phi (\alpha )$.
\end{defn}

In particular, since $i_{\partial} \geq -1$ and each
positive intersection contributes $2$ to the sum of intersection numbers, there can be 
no interior intersections with positive sign.
So we can reinterpret the definition:

\begin{lem}
An arc, $\alpha$, is sobering if and only if, after minimizing geometric 
intersections, $i_{\partial} \leq 0$, 
there are no positive (internal) intersections of $\alpha$ with $\phi 
(\alpha )$, and 
$\alpha$ isn't isotopic to $\phi(\alpha)$.
\end{lem}

There is a good reason to be interested in such arcs:

\begin{thm}
\label{sobering}
If there is a sobering arc $\alpha \subset \Sigma$ for $\phi$, then the open book $(\Sigma , \phi )$ 
is overtwisted.
\end{thm}

To prove this theorem we will construct a surface with Legendrian boundary which violates the 
Bennequin inequality, theorem \ref{bennineq}. Since this inequality is a necessary condition for 
tightness, 
we'll conclude that $\xi_{\Sigma, \phi}$ is overtwisted. In doing so, however, we will not find an 
overtwisted disk -- such a disk must exist, but will, in general, not be nicely positioned with 
respect to the open book. 

We begin our construction by suspending $\alpha$ in the mapping torus: 
$\alpha \times I \subset \Sigma 
\times_{\phi} I$. Let $p \in \partial \Sigma$ be an endpoint of 
$\alpha$. Since $\phi$ fixes the 
boundary, $p \times I$ is a meridian around the binding. 
Let $D_{\alpha}$ be the surface made by gluing a meridional disc to 
$\alpha \times I \subset M_{\Sigma , \phi}$ along $p \times I$ for 
each $p \in \partial \alpha$
(this extends $\alpha \times I$ across the binding in 
$M_{\Sigma , \phi}$). The surface $D_{\alpha}$ is a disk which has 
$\partial D_{\alpha} \subset \Sigma$, and is embedded on its 
interior, but possibly immersed on its boundary.

\begin{figure}[ht!]\anchor{fig:resolution}
\centering
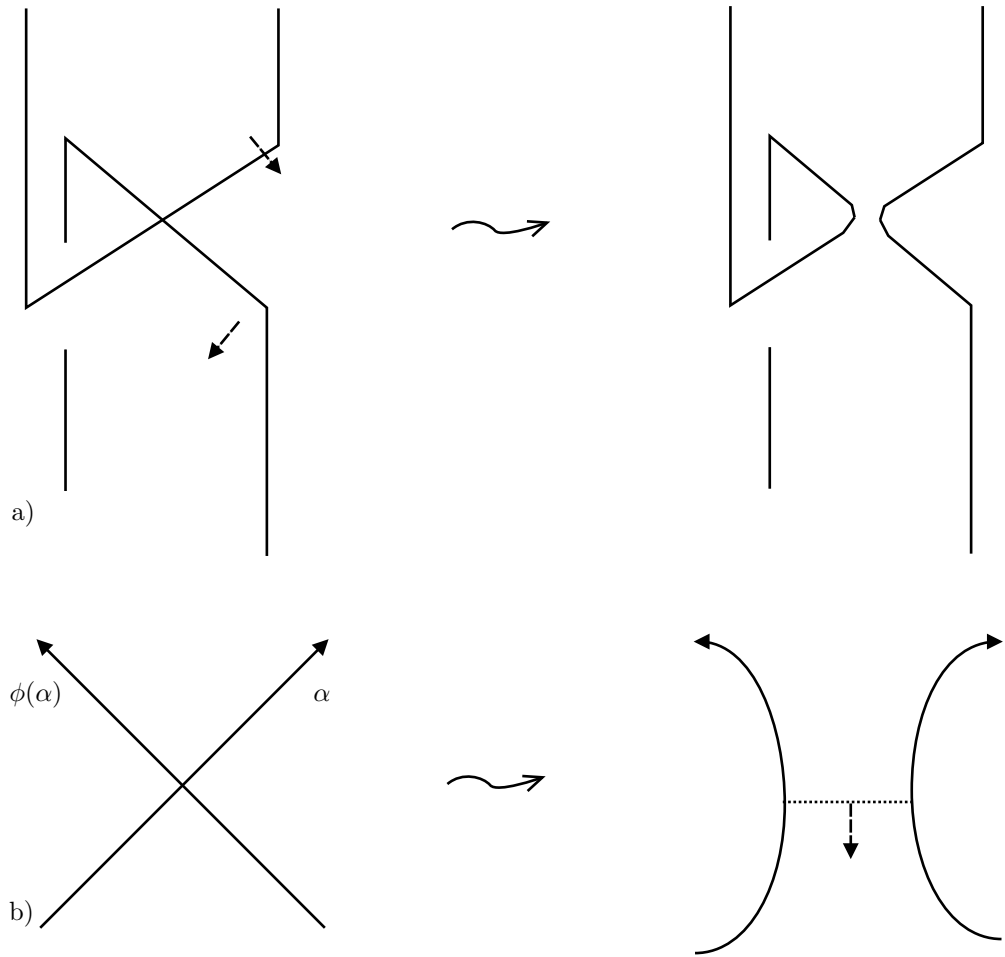
\caption{The local 
model for resolutions of immersed 
points on the boundary of a surface:\qua a) the 
surface, with dashed arrows indicating orientations, b) the resolution 
of $\alpha \cap \phi (\alpha 
)$, 
drawn on $\Sigma$, for the case of $D_{\alpha}$ (the dotted line 
represents $Int( D_{\alpha} ) \cap 
\Sigma$)}
\label{fig:resolution}
\end{figure}

We isotope $\phi (\alpha )$, relative to the boundary, to have minimal 
and transverse intersection with $\alpha$.
This induces an isotopy of $D_{\alpha}$, which intersects itself at transverse double points on its 
boundary. To 
proceed we need to remove the double points of $D_{\alpha}$. However, if we simply separate the 
edges of the disk meeting at an intersection, the boundary of the disk won't stay on $\Sigma$, 
eliminating any control we might have had over the framing. Instead we will opt to resolve the 
double points by decreasing the Euler characteristic of the surface:

\begin{defn}
Let $F$ be an orientable surface which is embedded on $Int(F)$ and has isolated double points on 
$\partial F$. The {\it resolution } of $F$ is constructed by thickening each double point into a 
half-twisted band in such a way that the resulting surface remains 
oriented. This is done according to the following local model (see 
\figref{fig:resolution}). In $\rr^3$ coordinates, in a small 
neighborhood of a double point, $F$ looks like $\{ (x,y,z) \in \rr^3 | z 
\geq 0, y = 0 \} \cup \{ (x,y,z) \in \rr^3 | z \leq 0, x = 0 \}$, 
oriented by $\hat{x}$ and $\pm \hat{y}$, respectively. We resolve by 
gluing in the triangles $ \{ (x,y,z) \in \rr^3 | z = 0, x \geq 0, 0 
\geq \pm y \geq x-1 \}$ and $ \{ (x,y,z) \in \rr^3 | z = 0, x \leq 0, 0
\leq \pm y \leq x+1 \}$ (that is, a pair of triangles depending on the 
intersection sign). Note that the boundary of both $F$ and its 
resolution are on the plane $z=0$.
\end{defn}

\begin{figure}[ht!]\anchor{fig:bdrysmooth}
\centering
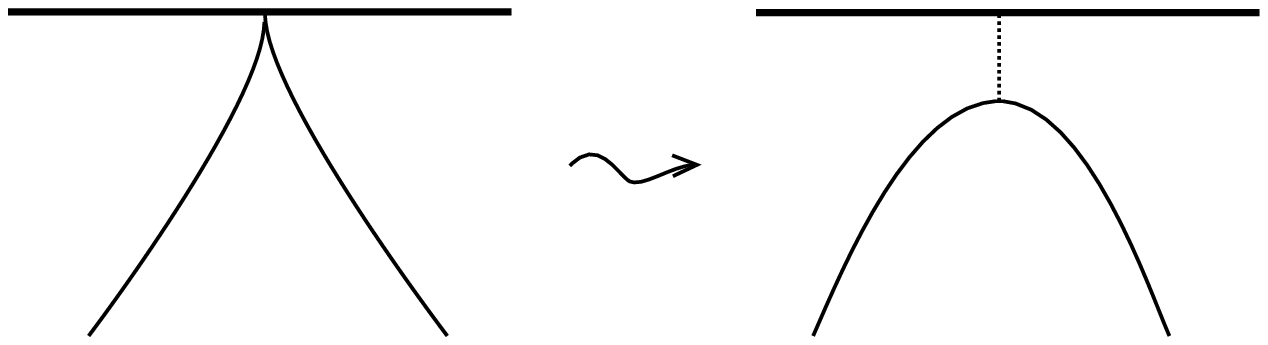
\caption{Smoothing the 
boundary of the suspension 
disk, $D_{\alpha}$:\qua The dashed line is $D_{\alpha} 
\cap \Sigma$.}
\label{fig:bdrysmooth}
\end{figure}

Let $S$ be the resolution of $D_{\alpha}$, with additionally 
$\partial D_{\alpha}$ pushed into 
$Int(\Sigma )$ at the points $\partial \alpha$ (see \figref{fig:bdrysmooth}). 
\figref{fig:resolution}b) shows, on $\Sigma$, how the boundary $\partial D_{\alpha}$ changes into 
$\partial S \subset \Sigma$ at a resolution.

\begin{lem}
$\chi (S) = 1 - i_{geom}(\alpha , \phi (\alpha ))$.
\end{lem}
\begin{proof}
The construction of $S$ starts with the disk $D_{\alpha}$ which has $i_{geom}(\alpha , \phi (\alpha 
))$ double points on its 
boundary. At each double point we glue a band from the disk to itself. The counting follows since a 
disk has $\chi = 1$, and adding a band reduces the Euler
characteristic by $1$.
\end{proof}

\begin{figure}[ht!]\anchor{fig:bdryframing}
\centering
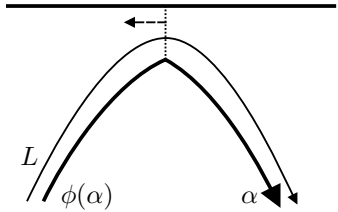
\caption{Framing near boundary: an example of computing the 
intersection of $L$ with $S$, on $\Sigma$. 
Here a positive boundary intersection at the start of $\alpha$ contributes $-1$ to the framing. The 
dotted line is $S \cap \Sigma$, while the dashed arrow is the orientation of $S$.}
\label{fig:bdryframing}
\end{figure}

\begin{lem}
The difference in framings of $\partial S$ with respect to $S$ and $\Sigma$ is: $Fr(\partial 
S; S,\Sigma ) = 
- i_{alg}(\alpha , \phi (\alpha )) - 
i_{\partial}(\alpha , \phi (\alpha ))$.
\end{lem}
\begin{proof}
The relative framing can be computed as the oriented intersection number of $S$ with a push-off of 
$\partial S$ 
along $\Sigma$. To choose the push-off first choose an orientation for $S$, and 
let
$s$ be a 
tangent vector to $\partial S$ which gives the boundary orientation. Choose 
vector $t \in T\Sigma$ so that $(s, t)$ provide an oriented basis agreeing with the orientation of 
$\Sigma$. Let $L$ be the push-off of $\partial S$ in the $t$ direction.

Now, $L \subset \Sigma$ can intersect $S$ only on $Int(S) \cap \Sigma$. This set consists of arcs 
where intersections were resolved in the construction of $S$ (including the boundary intersections, 
which are smoothed, as in \figref{fig:bdrysmooth}). 
For each of the internal intersections of $\alpha$ with $\phi(\alpha )$, $L$ intersects $S$ 
once with sign opposite that of the intersection, as shown in \figref{fig:framing}.

\begin{figure}[ht!]\anchor{fig:framing}
\centering
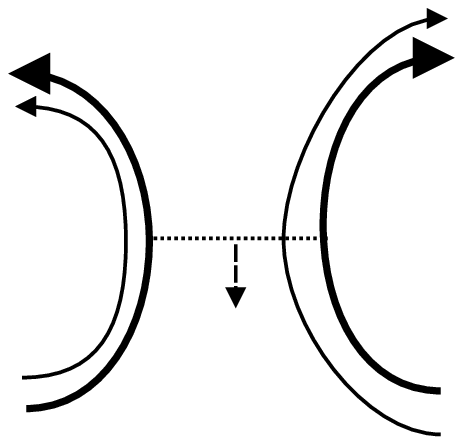
\caption{Framing at the resolution of a 
positive intersection:\qua The dashed line is $S \cap 
\Sigma$, while the arrow is the orientation of $S$. $L$ is a push-off of $\partial S$.}
\label{fig:framing}
\end{figure}

It remains to account for the boundary intersections. Since $\phi$ fixes $\partial \Sigma$ (as 
for any open book), there are exactly two intersections
points in $\partial \alpha \cap \partial \phi (\alpha )
= \partial \alpha$ (so $i_{\partial} (\alpha , \phi (\alpha )) \in \{ +1, -1, 0 \} $). 
There are four possibilities for the intersection of $L$ with $S$ at these resolutions, according 
to the signs of the intersections at the beginning and end of $\alpha$. 
One can compute each of these possibilities, as in the example of 
\figref{fig:bdryframing}, the result is shown in table \ref{table:bdrynumbers}.
There will be one intersection near the 
boundary, with sign $\mp 1$, if $i_{\partial}(\alpha , \phi (\alpha )) = \pm 1$, and 
either two with opposite signs, or none, if $i_{\partial}(\alpha , \phi (\alpha )) = 0$. Therefore the 
boundary intersections 
contribute $- i_{\partial}(\alpha , \phi (\alpha ))$ to the framing, $Fr(S,\Sigma )$.

Finally, in the above analysis we chose an orientation on $\alpha$  
(which then orients $\phi (\alpha )$, $D_{\alpha}$ and $S$). Reversing 
this orientation is equivalent to reversing the orientation of 
$\partial S$, and the effect on the framing computation is to take the push-off of $\partial S$ along 
$\Sigma$ but to the 
other side of $S$. However, the relative framing $Fr(S, \Sigma )$ is well-defined independent of which 
push-off is used, so our result doesn't depend on the chosen orientation. (The independence can also 
be seen directly by repeating the calculation using the other push-off.)
\end{proof}

\begin{table}[ht!]
\begin{center}
\begin{tabular}{|l|l|l|}
\cline{1-3}
& start $\alpha$ & end $\alpha$ \\
\cline{1-3}
pos. & -1 & 0 \\
neg. & 0 & +1 \\
\cline{1-3}
\end{tabular}
\caption[Framing intersections contributed by a boundary intersection.]{The framing intersections 
contributed by a boundary intersection of the given sign at the
beginning or end of $\alpha$. Computations were made as in \figref{fig:bdryframing}.}
\label{table:bdrynumbers}
\end{center}
\end{table}

We are now ready to prove the main theorem of this section.

\begin{proof}[Proof of Theorem \ref{sobering}]

As we will later wish to LeRP $\partial S$, we must first rule out any 
components of $\partial S$ which bound 
disks on $\Sigma$. 

Assume there is a disk $D \subset \Sigma$ such that $\partial D \subset \partial S$. The intersection 
signs at the corners of a disk (where resolutions have occurred on $\partial S$) must alternate: see 
\figref{fig:diskalt}.
However, by the sobering condition, there is at most one positive intersection, which occurs at the 
boundary. This implies that the only possible disk is a bi-gon with a corner on the boundary, but in 
this case either the intersection $\alpha \cap \phi (\alpha )$ isn't minimal, or 
the other corner of the bi-gon is also on the boundary. The former contradicts the minimal set-up 
we've arranged, 
the latter contradicts the definition ($\alpha$ isn't isotopic to $\phi (\alpha )$), so 
there 
can be no such disk. 

\begin{figure}[ht!]\anchor{fig:diskalt}
\centering
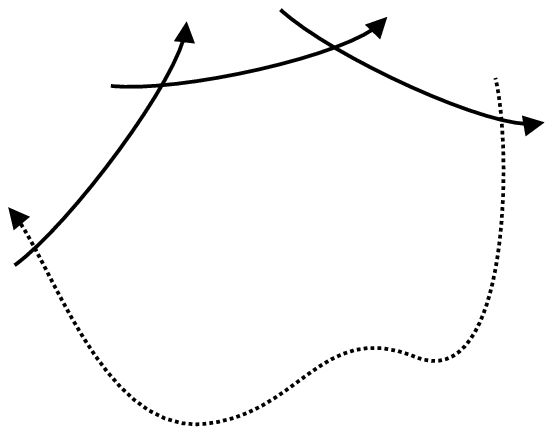
\caption{If disk $D \subset \Sigma$ was created by resolving $\alpha 
\cap \phi (\alpha )$, the signs at the 
intersections would have to alternate.}
\label{fig:diskalt}
\end{figure}

In order to use the machinery of convex surfaces, we must make the leaves convex. We can do this by an isotopy of the contact structure.

\begin{lem}
\label{convexleaf}
Let $(\Sigma ,\phi )$ be an open book, $\pi : M_{\Sigma ,\phi} \setminus
\partial \Sigma \rightarrow S^1$ the fibration of the complement to the binding.
The closure of $\pi^{-1}(0) \cup \pi^{-1}(0.5)$ 
is convex and divided by the binding, $\partial \Sigma = \partial (\pi^{-1}(0))$, for some compatible contact 
structure.
\end{lem}
\begin{proof}
Let $\tilde{\Sigma}$ be the closure of $\pi^{-1}(0) \cup \pi^{-1}(0.5)$ 
(oriented by $\pi^*(ds)$ on $\pi^{-1}(0)$ and 
$\pi^*(-ds)$ on $\pi^{-1}(0.5)$), we take a 
compatible contact structure built as in \cite{thurston-winkel} from a 
one-form 
$\lambda$ (see \cite{diss}). We may choose the homotopy extending 
$\lambda$ to the mapping torus 
in such a way that the characteristic foliation on $\tilde{\Sigma}$ 
looks like $ker(\lambda + \phi^*(\lambda))$ on 
$\pi^{-1}(0.5)$ (and like $ker(\lambda)$
on $\pi^{-1}(0)$). By a small perturbation of this homotopy (eg. by 
perturbing $\phi$) we can 
assure that the saddles on $\pi^{-1}(0)$ connect to the sources/sinks on 
$\pi^{-1}(0.5)$, and vice versa, hence 
the characteristic foliation is Morse-Smale. By the results of 
\cite{girouxconvex}, a surface with such a characteristic foliation is 
convex. The dividing curves 
separate the regions where the contact one-form is positive on a 
positive normal vector from those where it is 
negative. Since the contact structure is compatible, $\pi^{-1}(0)$ is 
positive while 
$\pi^{-1}(0.5)$ is negative (with the natural orientation: $\pi^*(ds)$, 
on each), hence 
the dividing curve is exactly the binding. 
\end{proof}

Because the compatible contact structure to an open book is unique up to isotopy, we can perform an isotopy 
to assure that our contact structure is that of the lemma.
Now, with convex pages, we want to apply the LeRP (lemma \ref{lerp}) to 
the curves $\partial S$, but we must insure that 
there is a component of the dividing set in each component of  
$\tilde{\Sigma} \setminus \partial S$. 
To 
achieve this, let $\Gamma$ contain a simple closed curve in each 
component of $\tilde{\Sigma} \setminus 
\partial S$. Using the folding lemma \ref{folding} along $\Gamma$ we 
perturb $\tilde{\Sigma}$, changing 
the dividing set to $\partial \Sigma \cup \Gamma \cup \Gamma'$ 
($\Gamma'$ a parallel copy of 
$\Gamma$). 
Now we can apply the LeRP, perturbing $\tilde{\Sigma}$ again to make the 
curves $\partial S$ 
Legendrian.

By the definition of convex surfaces there is a contact vector field $v$ 
normal to $\tilde{\Sigma }$.
Since $\partial S$ doesn't cross the dividing set, $v |_{\partial S} \not\in \xi_{\Sigma, \phi }$, 
so the framing of $S$ with
respect to $\xi_{\Sigma, \phi }$ (the twisting number)
may be computed by using a push-off of $\partial S$ in the $v$ 
direction. But $v$ is normal to 
$\tilde{\Sigma }$ so a push-off along $\tilde{\Sigma}$ is equivalent, 
for framing purposes, to the push 
off along $v$. Hence, the framing coming from $\xi_{\Sigma , \phi}$ is 
the same as that from $\Sigma$: $tw(\partial S, S) = Fr(\partial S; S,\Sigma )$.

We computed above that $Fr(S, \Sigma) =  - i_{alg}(\alpha , 
\phi (\alpha )) - i_{\partial}(\alpha , \phi (\alpha ))$. So,
\[ tw(\partial S, S) + \chi (S) \geq -i_{alg}(\alpha ,
  \phi (\alpha )) - i_{\partial}(\alpha , \phi (\alpha ) + 1 - i_{geom}(\alpha , \phi (\alpha )). \]
The sobering condition $-i_{alg}(\alpha ,
  \phi (\alpha )) - i_{\partial}(\alpha , \phi (\alpha ) - i_{geom}(\alpha , \phi (\alpha )) \geq 0$ 
then implies:
\[ tb(\partial S, S) + \chi (S) > 0, \]
violating the Bennequin inequality (theorem \ref{bennineq}), which must hold if $\xi_{\Sigma ,\phi}$ 
is 
tight. Hence, $(\Sigma ,\phi )$ is overtwisted.
\end{proof}

\begin{rem}
A word about the use of arcs instead of closed curves inside $\Sigma$. Our construction, above,
would fail for these curves for simple Euler characteristic reasons, or we may
rule out such considerations more directly: in the complement of the binding $\xi_{\Sigma, \phi}$ is
a perturbation of the taught foliation given by the leaves of the open book. From the results of
Eliashberg and Thurston \cite{confoliations} such a contact structure is symplectically fillable, and
hence tight. Thus, no surface entirely in the complement of the binding, such as the suspension of a 
closed curve on $\Sigma$, can violate the Bennequin inequality.
\end{rem}

We can slightly extend the class of monodromies for which sobering arcs are useful by a little 
cleverness:

\begin{thm}
If an arc $\alpha \subset \Sigma$ has:
\[ i_{alg}(\alpha , \phi (\alpha
)) + i_{\partial}(\alpha , \phi (\alpha )) + i_{geom}(\alpha , \phi (\alpha )) = -1, \]
 then each open 
book $(\Sigma, 
\phi^n)$, for $n>0$, is overtwisted.
\end{thm}
\begin{proof}
The open book $(\Sigma, \phi^n)$ can be built by gluing $n$-copies of $\Sigma 
\times I$ to each other 
by $\phi$ (including the bottom of the first to the top of the last) to get the mapping torus $\Sigma 
\times_{\phi^n} S^1$, then adding the binding as usual. In each $\Sigma \times I$ we have a disk 
$\alpha 
\times I$. The top of one disk intersects the bottom of the next as $\alpha \cap \phi (\alpha)$. If we 
resolve these intersections, as in the previous proofs, we get a surface with $\chi = n - n i_{geom} 
(\alpha, \phi (\alpha))$ and framing $n i_{alg}(\alpha, \phi (\alpha))$. Next we add in the binding 
and cap off the surface with disks $D^2 \times \partial \alpha$, to get a surface $S$ with $\chi (S) 
= 2 -n - n i_{geom}(\alpha, \phi (\alpha))$ and framing $n i_{\partial}(\alpha, \phi (\alpha)) + n 
i_{alg}(\alpha, \phi (\alpha))$. 

To finish the proof we proceed exactly as in the proof of theorem 
\ref{sobering}: rule out components of 
$\partial S$ which bound disks on a page, then fold and LeRP to make $\partial S$ Legendrian (we'll 
have to do this along several pages, but they won't interfere since folding and LeRPing can be done in 
small neighborhoods). Finally we get a surface with 
$tb(\partial S, S) + \chi (S) = -n (i_{alg}(\alpha ,
\phi (\alpha )) + i_{\partial}(\alpha , \phi (\alpha )) ) +2 -n$. Then under the condition 
$i_{alg}(\alpha , \phi (\alpha
)) + i_{\partial}(\alpha , \phi (\alpha )) + i_{geom}(\alpha , \phi (\alpha )) = -1$ we 
get:
$tb(\partial S, S) + \chi (S) = n+2-n =2 > 0$, so the open book is overtwisted.
\end{proof}

We can apply this theorem to get some interesting examples:

\begin{cor}
\label{boundneg}
Let $\Sigma_g$ be a surface with genus $g > 0$ and one boundary 
component, and let $\delta \subset 
\Sigma_g$ be parallel to the boundary. The open book $(\Sigma_g, D_{\delta}^{-n})$ is overtwisted 
($n>0$).
\end{cor}
\begin{proof}
Let $\{ a_i, b_i \} _{i=1,\dots , g}$ be simple closed curves on $\Sigma_g$ such that $a_i$ 
intersects $b_i$ in one point, $b_i$
intersects $a_{i+1}$ in one point, and the curves are otherwise disjoint.  
One can show, as in \cite{piergalliniloi}, that $D_{\delta} = (\prod_i{D_{a_i}D_{b_i}})^{4g+2}$.

Then $D_{\delta}^{-n} = (D_{b_g}^{-}D_{a_g}^{-}\cdots D_{b_1}^{-}D_{a_1}^{-})^{n(4g+2)}$. An arc which intersects 
only $b_g$ will be sobering for $\phi = D_{b_g}^{-}D_{a_g}^{-}\cdots D_{b_1}^{-}D_{a_1}^{-}$.
\end{proof}

On the other hand, since the open books $(\Sigma_g, D_{\delta}^{n})$ are positive, 
they are Stein fillable, and hence tight. (The connection between fillable contact structures and 
positive monodromy was found by many people, see \cite{steinpos} for a review.)

We return now to the Poincare
homology-sphere with reversed orientation, discussed in the introduction.
This manifold arises as $+1$-surgery on the trefoil $T^+ \subset S^3$, 
so it supports an open book which comes from the natural 
fibering of $S^3 \setminus T^+$ by adding a boundary parallel negative Dehn twist.
 (To see this, push the surgery curve onto a 
page, then verify directly that an arc passing the surgery experiences a 
negative Dehn twist). 
Explicitly: let $\Sigma$ be the 
punctured torus, $\delta$ a boundary parallel curve, and 
$a,b$ curves which intersect in a point. Then the open book $(\Sigma, 
D^{-1}_{\delta}D^{+1}_{b}D^{+1}_{a} )$ represents the Poincare 
homology-sphere with reversed
orientation.
The technique of the corollary shows that the monodromy can be written $\phi = 
(D_b^- D_a^-)^6(D_b^+ D_a^+) = (D_b^- D_a^-)^5$, and that this is overtwisted.

\section{A special case}
\index{Detecting Hopf Bands}
\label{specialcase}

\begin{figure}[ht!]\anchor{fig:hopf}
\centering
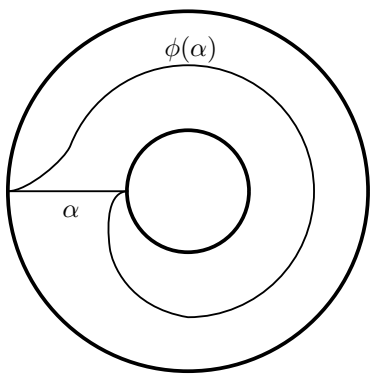
\caption{A negative Hopf-band, $H^-$, with the transverse 
arc, $\alpha$, and its image}
\label{fig:hopf}
\end{figure}

Let us consider the simplest open books. In a disk all arcs are boundary parallel, so we turn our 
attention to annuli, and especially to the Hopf-bands, $H^{\pm}$, which 
are annuli with monodromy a single 
right or left-handed Dehn twist along, $c$, the center of the annulus. 
If we 
let $\alpha$ be an arc which 
crosses the 
annulus, we see from \figref{fig:hopf} that $i_{geom} ( \alpha , 
D^\pm_c (\alpha )) = 
i_{alg} ( \alpha , D^\pm_c (\alpha )) = 0$, while $i_{\partial} (\alpha 
, D^\pm_c (\alpha )) = \pm 1$ for 
$H^{\pm}$, respectively. Since $D^\pm_c (\alpha )$ isn't isotopic to 
$\alpha$, $\alpha$ is 
sobering for $H^-$ (but not for $H^+$!).

\begin{lem}
\label{neghopfsum}
The Murasugi sum $(\Sigma , \phi ) * H^-$, for any open book $(\Sigma , 
\phi )$, along any attaching curve, is overtwisted.
\end{lem}
\begin{proof}

Take an arc across $H^-$, outside the attaching region. It 
is unmoved by (the extension of) $\phi$, so will have the same intersection properties 
relative to $D^-_c \circ \phi$ as 
the arc $\alpha$, above, had relative to $D^-_c$. So this arc will be 
sobering in the sum.
\end{proof}

How special are the intersection properties of the arcs we've just considered? It turns out that they 
uniquely specify the geometric situation of having a Hopf-band Murasugi summand.

\begin{thm}
\label{hopfband}
If $i_{geom} ( \alpha , \phi (\alpha )) = 0$ and $i_{\partial} (\alpha , \phi (\alpha )) = \pm 1$, 
then $\alpha$ is the transverse arc to a positive (respectively negative) Hopf-band. That is $(\Sigma 
, \phi ) = (\Sigma ', \phi ') * H^{\pm}$ where $\alpha$ is parallel to the attaching curve in 
$H^{\pm}$.
\end{thm}
\begin{proof}
Let $c \subset \Sigma$ be the curve which comes from smoothing $\alpha \cup \phi (\alpha 
)$. Since $i_{geom} ( \alpha , \phi (\alpha )) = 0$, $c$ will be a simple closed curve. Note that 
$\phi (\alpha ) = D^{\pm}_c (\alpha )$ up to isotopy, where the sign of the Dehn-twist 
agrees with $i_{\partial} (\alpha , \phi (\alpha ))$.

Let $\phi ' = D^{\mp}_c \circ \phi$ on the surface $\Sigma ' = \Sigma \setminus \alpha$. This makes 
sense since $\phi '(\alpha ) = D^{\mp}_c \circ \phi (\alpha ) =  D^{\mp}_c \circ D^{\pm}_c (\alpha ) = 
\alpha$ (up to isotopy, so for some representative of the isotopy class of $\phi'$).

If we plumb a Hopf-band onto $(\Sigma ', \phi ')$ along $c \setminus \alpha \subset \Sigma '$, we get 
the map $D^{\pm}_c \circ \phi ' = \phi$ on $\Sigma$. Thus $(\Sigma , \phi ) = (\Sigma ', \phi ') * 
H^{\pm}$, and $\alpha$ crosses the one-handle, as required.
\end{proof}

Of course if  $(\Sigma , \phi ) = (\Sigma ', \phi ') * H^{\pm}$, then there is such an arc $\alpha$, 
described above, 
so we in fact have a necessary and sufficient condition for the existence of a Hopf-band 
summand in an open 
book. This criterion can be used to show that certain fibered links are not Hopf plumbings, and it is 
often helpful in finding stabilizations in otherwise difficult examples, as in section \ref{sympl}.

\section{A criterion for overtwisting}
\label{criterion}

\begin{thm}                                                              
\label{sobcrit}
An open book is overtwisted if and only if it is $H^+$-stably equivalent 
to an open book with a sobering arc.
\end{thm}  
\begin{proof}                                                  
{\it If:}\qua The final open book of the equivalence has a sobering arc so 
its compatible contact 
structure is overtwisted. Neither summing nor de-summing $H^+$ changes 
the contact structure, so the original open book is overtwisted.
                                             
{\it Only If:}\qua 
For any open book $B$ the homotopy classes 
of plane fields on $M_{B}$ have a natural $\zz$--action,
and we arrive in the homotopy class of $\xi_{B*H^-}$ 
by acting on the homotopy class of $\xi_B$ with $+1$ (see \cite{girouxgoodman}).

Now let $(\Sigma, \phi)$ be an overtwisted open book, let $\xi$ be a 
contact structure in the homotopy class of 
$\xi_{(\Sigma, \phi)}$ acted on by $-1$ (there is a contact structure in 
any homotopy class), and let $B$ be an open book compatible 
with $\xi$. The contact structure $\xi_{B * H^-}$ then is homotopic to 
$\xi_{(\Sigma, \phi)}$, and both are overtwisted, hence they are 
isotopic (lemma \ref{otclass}). So, $B * H^-$ is $H^+$-stably equivalent 
to $(\Sigma, \phi)$ by theorem \ref{stabeq}, 
and $B * H^-$ has a sobering arc (across the Hopf band, see lemma \ref{neghopfsum}).
\end{proof}

\section{Boundaries of symplectic configurations}
\label{sympl}
\index{Boundaries of Symplectic Configurations}

Following \cite{gay} we define a {\it symplectic configuration graph} 
to be a labeled graph $G$ with 
no edges from a vertex to itself and with each vertex $v_i$ labeled 
with a triple $(g_i,m_i,a_i)$, 
where $g_i \in \{ 0,1,2,\dots \} $, $m_i \in \zz$ and $a_i \in 
(0,\infty )$. ($a_i$ doesn't enter 
into our purely 
three-dimensional concerns, so we'll often omit it.) Let $d_i$ denote 
the degree of vertex $v_i$. A 
configuration graph is called {\it positive} if $m_i + 
d_i > 0$ for every vertex $v_i$.

\begin{figure}[ht!]\anchor{fig:confgraph}
\centering
\includegraphics[scale=0.8]{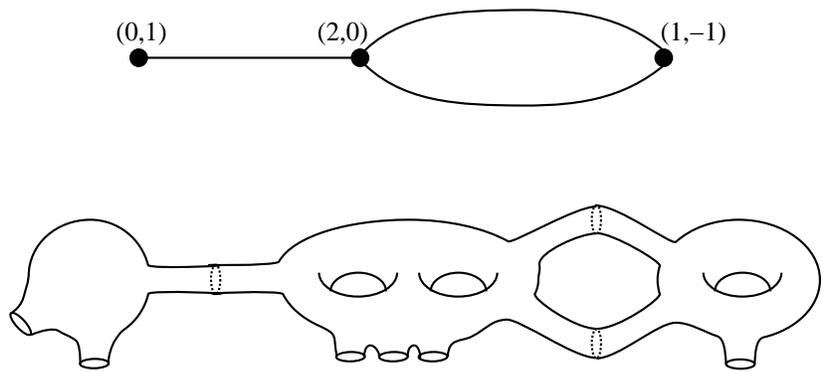}
\caption{A positive configuration graph (the $a_i$ have been 
suppressed), and the surface 
$\Sigma (G)$:\qua The dotted circles are the $e_{ij}$.}
\label{fig:confgraph}
\end{figure}

Given a positive configuration graph $G$ we define an open book \\
$(\Sigma (G), \phi (G))$ as follows: 
For each vertex $v_i$ let $F_i$ be a surface of genus $g_i$ with $m_i 
+ d_i$ boundary components. 
$\Sigma (G)$ is the surface obtained by connect sum of the $F_i$, with 
one connect sum between $F_i$ 
and $F_j$ for each edge connecting $v_i$ to $v_j$. See \figref{fig:confgraph} for an example. 
For each 
edge in G there is a circle $e_{ij}$ in $\Sigma (G)$. Let 
$\sigma (G) = \prod D^+_{e_{ij}}$, a right-handed Dehn twist at each connect sum. Let $\delta (G)$ 
be the product of one right-handed Dehn twist around each circle of $\partial \Sigma (G)$. Finally, 
$\phi (G) = \sigma (G)^{-1} \circ \delta (G)$.

In \cite{gay} Gay shows that the open books $(\Sigma (G), \phi (G))$ arise as the concave boundary of 
certain symplectic manifolds. In that paper he asks the question:

\begin{quest}(Gay \cite{gay})\qua
Are there any positive configuration graphs $G$ for which we can show
that $(\Sigma (G), \phi (G))$ is overtwisted (and hence conclude that a symplectic configuration with 
graph $G$ cannot embed in a closed symplectic 4-manifold)?
\end{quest}

We can use the techniques of the previous sections to give such examples.

First, note that if $l$ is a properly embedded arc in $\Sigma (G)$ which crosses from one $F_i$ to an 
adjacent one, so 
crosses one negative Dehn twist and two positive, then $i_{geom} ( l , \phi (l)) = 0$ and 
$i_{\partial} (l , \phi (l)) = +1$ (see \figref{fig:lhopf}). From theorem \ref{hopfband} we know 
that there is an $H^+$ summand. To remove it we'll first re-write the monodromy using the lantern 
relation (this is a standard fact of mapping class groups, see \cite{mapclass}):

\begin{figure}[ht!]\anchor{fig:lhopf}
\centering
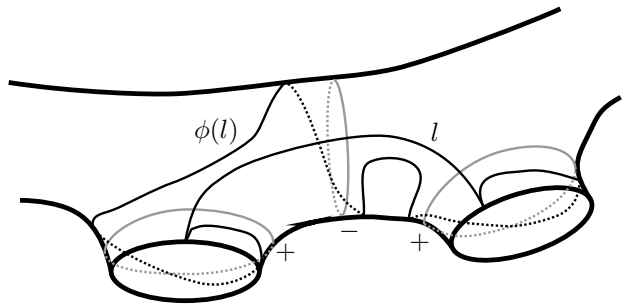
\caption{The arc $l \subset \Sigma (G)$, its image under $\phi (G)$, 
and the twisting curves it 
crosses}
\label{fig:lhopf}
\end{figure}

\begin{figure}[ht!]\anchor{fig:deplumb}
\centering
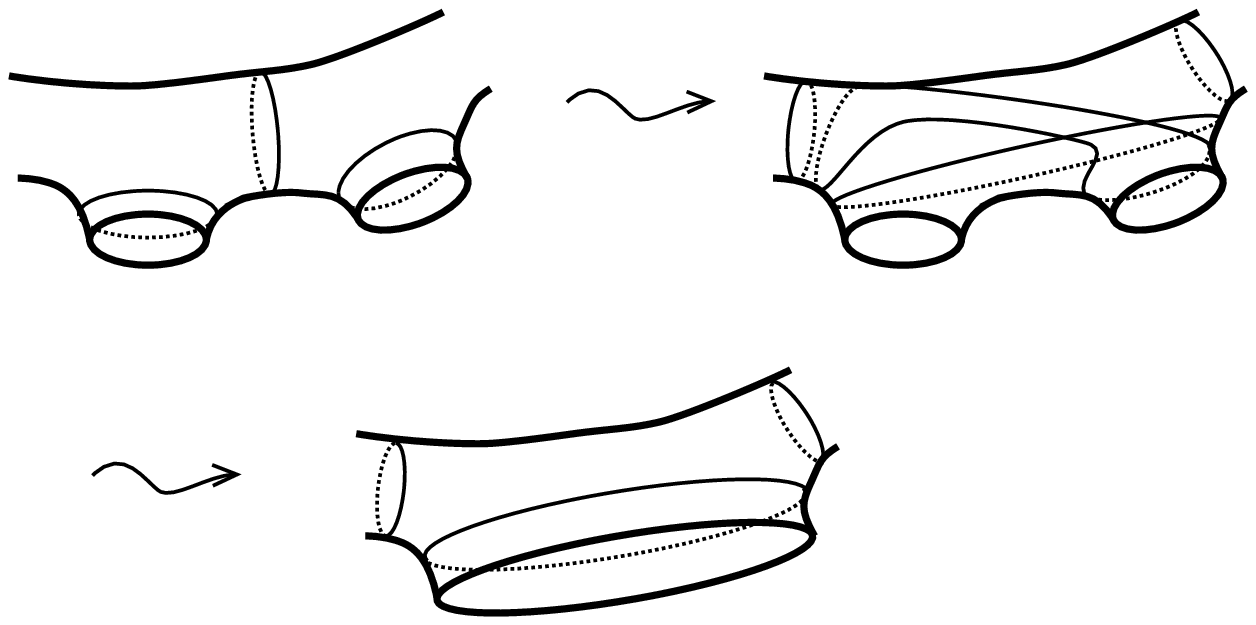
\caption{De-plumbing an $H^+$ from $(\Sigma (G), 
\phi (G))$. The first step is an application of the
lantern lemma to go from $D^+_{a} \circ D^+_{b} \circ D^-_{\alpha}$ to
$D^-_{c} \circ D^-_{d} \circ D^+_{\beta} \circ D^+_{\gamma}$,
the second is the removal of a Hopf-band along $\gamma$.}
\label{fig:deplumb}
\end{figure}

\begin{lem}
\label{lantern}
The following relation holds in the mapping class group of a four-times punctured sphere (fixing 
the boundary point-wise), referring to \figref{fig:deplumb}(a) and (b):
\[ D^+_{a} \circ D^+_{b} \circ D^+_{c} \circ D^+_{d} = D^+_{\alpha} \circ D^+_{\beta} \circ 
D^+_{\gamma}. \]
\end{lem}

Applying this relation to the region of \figref{fig:lhopf} around arc $l$, we can rewrite 
the monodromy as in \figref{fig:deplumb} (a) and (b): 
$D^+_{a} \circ D^+_{b} \circ D^-_{\alpha} = D^-_{c} \circ D^-_{d} \circ D^+_{\beta} \circ 
D^+_{\gamma}$ (Dehn twists about $a$, $b$, $c$, $d$ commute with each other and with $\alpha$, 
$\beta$, and $\gamma$, since these pairs don't intersect). 
It is then clear, since $l$ crosses only the final Dehn twist $D^+_{\gamma}$, how to 
remove a 
Hopf-band: remove $D^+_{\gamma}$ and cut along $l$, \figref{fig:deplumb}(c). 
We arrive at a surface with one fewer boundary component, and $D^-_{c} \circ D^-_{d} 
\circ D^+_{\beta}$ replacing $D^+_{a} \circ D^+_{b} \circ D^-_{\alpha}$ in the monodromy. We'll use 
this 
move repeatedly to prove the 
next lemma. Note that this move doesn't necessarily return an open book $(\Sigma (G'), \phi 
(G') )$, since the resulting open book may have several negative Dehn twists around the same curve.

\begin{thm}
\label{posconfex}
Let $G$ be a positive configuration graph with a vertex 
$v_1$ such that $d_1 = 1$ and $g_1 = 0$, but not a graph with:
\begin{enumerate} \nonumber
 \item one edge, $m_1 = 0$,
 \item one edge, $g_2 = m_2 = 0$,
 \item one edge, $g_2 = 0$, $m_2 = m_1 = 1$, or,
 \item two edges, $g_3 = 0$, $m_1 = m_3 = 0$.
\end{enumerate} 
Then $(\Sigma (G), \phi (G))$ is overtwisted.
\end{thm}
\begin{proof}
First we apply the move illustrated in \figref{fig:deplumb} 
repeatedly, removing positive 
Hopf-bands until all of the 
boundary components of $F_1$ are gone, leaving a boundary component 
with one positive Dehn twist and 
$m_1+1$ negative twists. These, of course, are equivalent to $m_1$ 
negative Dehn twists. We are left with a surface, $G'$, in which $F_1$ 
and $F_2$ have been replaced with a surface of genus $g_2$, and $m_2 + 
d_2$ boundary components, which has $m_1$ negative Dehn twists about one 
of these punctures (the remaining monodromy twists are unchanged).
See \figref{fig:deplumblots}.

\begin{figure}[ht!]\anchor{fig:deplumblots}
\centering
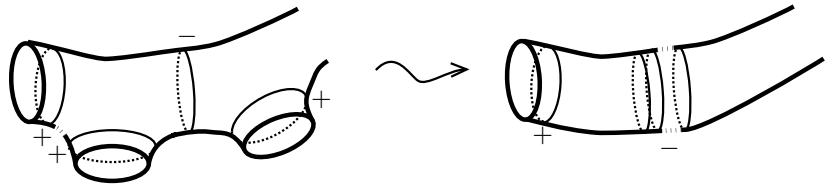
\caption{Removing almost all boundary 
components by applying Figure \ref{fig:deplumb} to remove $H^+$}
\label{fig:deplumblots}
\end{figure}

The surface $G'$ will have $\sum_{i>1} m_i + d_i$ boundary components and, 
since $m_i + d_i > 0$, will have a single boundary component only 
when $e=1$ and $m_2 = 0$ (where $e$ is the number of edges). In this case we have a surface with one 
puncture and monodromy consisting of $m_1$ negative twists about the 
boundary; by corollary \ref{boundneg} the open book will be 
overtwisted if $g_2 > 0$ and $m_1 >0$. The 'bad' cases of $e=1$, 
$m_2=0$, and either $g_2 = 0$ or $m_1 = 0$ are excluded by hypothesis.

Now if $G'$ has several boundary components, let $\alpha$ be an arc from 
the boundary component with 
$m_1$ negative Dehn twists to 
another boundary component, chosen so that it crosses as many of 
the $e_{ij}$ (non-boundary parallel negative twists of $\phi(G)$) as possible, but crosses each at most once. The arc $\alpha$ will cross the $m_1$ negative Dehn 
twists, $k$ additional negative twists (one for each $e_{ij}$ crossed), and finally a positive Dehn 
twist (see \figref{fig:soberex}). (Note that if $k=0$ then $e=1$, and if $k=1$ then $e=2$: if there were more edges then $\alpha$ could be extended.)

%(Note that $k \geq e-1$ since each possible $\alpha$ gives a path in the graph which starts at $v_2$ (the neighbor of $v_1$) and crosses one negative twist for each edge of the path.) 

Before any isotopy there are only negative internal intersections of 
$\alpha$ with $\phi(\alpha)$, thus after an isotopy which minimizes 
intersections there will still be only 
negative internal intersections. The boundary intersection number will 
initially be $0$ (since $\alpha$ crosses only the one positive twist as in \figref{fig:soberex}), and 
can only be reduced as internal intersections are minimized (since all 
internal intersections are negative). 
%The arc $\alpha$ crosses at least one negative twist except when $e=1$ and $m_1=0$ (see \figref{fig:soberex}). 
So, $\alpha$ will be sobering if it is not isotopic 
to $\phi(\alpha)$.

If $m_1 + k > 1$, that is $\alpha$ crosses more than one negative twist, 
the algebraic intersection number (including boundary 
intersections) is less than $0$. This is an isotopy invariant, so 
$\alpha$ is not isotopic to $\phi(\alpha)$.

We can rule out the case $m_1 + k = 0$ since we must have $e=1$ and $m_1=0$. 
What if $m_1 + k = 1$?  Then $\alpha$ crosses one negative and one positive Dehn twist, call them $D_{\beta}^+$ 
and $D_{\gamma}^-$. Up to homotopy (relative to the boundary) 
$[\phi(\alpha) ] = [\alpha] + [\beta ] +[\gamma ]$ (as $\beta$ and $\gamma$ can't cross). If $\alpha$ is isotopic to $\phi(\alpha)$ then 
$[\beta ]  = [\gamma ]$, up to sign.

There are two cases, first consider $m_1 = 1$ and $k=0$. In this case $e=1$, so 
both $\beta$ and $\gamma$ are boundary parallel -- in order for them to 
be homotopic $G'$ must be a cylinder, so $g_2=0$ and $m_2=m_1=1$ (a case 
excluded by hypothesis).

In the other case $m_1=0$ and $e=2$, now $\beta$ is again boundary parallel 
(near the boundary of $F_3$), but $\gamma$ is the curve separating $F_3$ 
from $F_2$. Their homotopy equivalence implies that $F_3$ is a 
cylinder, so $g_3=0=m_3$ (another excluded case).

\begin{figure}[ht!]\anchor{fig:soberex}
\centering
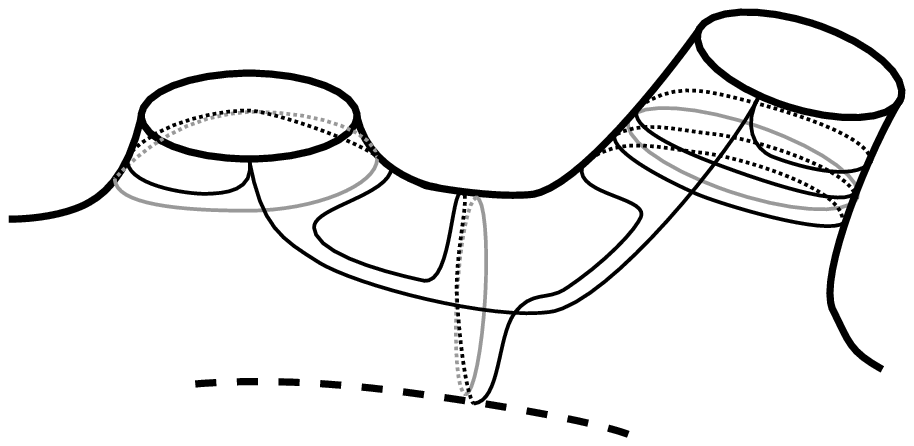
\caption{Schematic of the sobering curve}
\label{fig:soberex}
\end{figure}

We have shown that $\alpha$ will be sobering in the allowed cases. By 
theorem \ref{sobering}, the 
presence of a sobering arc implies that the open book is overtwisted. 
When stabilizing, to get back the open book $(\Sigma (G), \phi (G))$, the 
compatible contact structure remains overtwisted (in fact, the same up 
to isotopy), so $(\Sigma (G), \phi (G))$ is overtwisted.
\end{proof}

\begin{rem}
Since they are overtwisted these contact manifolds can't be the convex boundary of a symplectic four-manifold.
Corollary 2.1 of \cite{gay} also shows this for some (but not all) of these graphs. The techniques 
used there are four-dimensional -- a clever application of the adjunction inequality.
\end{rem}

\Addresses\recd

\end{document}